\documentclass[twoside,british,english]{elsarticle}
\usepackage[T1]{fontenc}
\usepackage[latin9]{inputenc}
\usepackage{geometry}
\usepackage{color}
\usepackage{textcomp}
\usepackage{amsmath}
\usepackage{amsthm}
\usepackage{amssymb}
\PassOptionsToPackage{normalem}{ulem}
\usepackage{ulem}

\journal{Journal of Mathematical Analysis and Applications}

\theoremstyle{plain}
\newtheorem{thm}{\protect\theoremname}
\theoremstyle{plain}
\newtheorem{prop}[thm]{\protect\propositionname}
\theoremstyle{definition}
\newtheorem{example}[thm]{\protect\examplename}
\newtheorem{rem}{Remark}
\newtheorem{cor}{Corollary}

\usepackage{textcomp}

\usepackage{babel}
\addto\captionsbritish{\renewcommand{\examplename}{Example}}
  \addto\captionsbritish{}
  \addto\captionsbritish{\renewcommand{\propositionname}{Proposition}}
  \addto\captionsbritish{}
  \addto\captionsbritish{\renewcommand{\theoremname}{Theorem}}
  \addto\captionsenglish{\renewcommand{\examplename}{Example}}
  \addto\captionsenglish{}
  \addto\captionsenglish{\renewcommand{\propositionname}{Proposition}}
  \addto\captionsenglish{}
  \addto\captionsenglish{\renewcommand{\theoremname}{Theorem}}
  \providecommand{\examplename}{Example}
  
  \providecommand{\propositionname}{Proposition}
  
\providecommand{\theoremname}{Theorem}

\makeatother

\begin{document}
\begin{frontmatter}
\selectlanguage{english}%

\title{ Isometric immersions and differential equations which describe pseudospherical surfaces}

\author[add1]{Diego Catalano Ferraioli\corref{cor1}}

\ead{diego.catalano@ufba.br}

\author[add2]{Tarc\'isio Castro Silva}

\ead{tarcisio@mat.unb.br}

\author[add3]{Keti Tenenblat}

\ead{k.tenenblat@mat.unb.br}

\cortext[cor1]{Corresponding author at: Instituto de Matem\'atica e Estat\'istica- Universidade Federal
da Bahia Campus de Ondina, Av. Adhemar de Barros, S/N, Ondina - CEP
40.170.110 - Salvador, BA - Brazil }

\address[add1]{Instituto de Matem\'atica e Estat\'istica- Universidade Federal da Bahia, Campus de
Ondina, Av. Adhemar de Barros, S/N, Ondina - CEP 40.170.110 - Salvador, BA - Brazil}

\address[add2]{Department of Mathematics, Universidade de Brasilia, Brasilia DF 70910-900, Brazil}

\address[add3]{Department of Mathematics, Universidade de Brasilia, Brasilia DF 70910-900, Brazil}

\begin{abstract}
We provide families of second order non-linear partial
differential equations describing pseudospherical surfaces (\textbf{pss} equations)
with the property of having local isometric immersions
in $\mathbb{E}^{3}$ with principal curvatures depending on the  finite order jets of solutions of the differential equations. These equations occupy a particularly special place amongst 
\textbf{pss} equations since a series of recent papers on several classes of \textbf{pss} equations seemed to suggest that only the sine-Gordon equation had the above property. Explicit examples are given, which include the short pulse equation and some generalizations.
 
\end{abstract}


\end{frontmatter}

\section{Introduction}
A differential equation $\mathcal{E}$, for $z=z(x,t)$, is said to
be a \textit{\textbf{pss} equation} (or \textit{to describe pseudospherical
surfaces}) if it is equivalent to the structure equations of a surface
with constant Gaussian curvature equal to $-1$, i.e., 
\begin{equation}\label{struct} 
d\omega_{1}=\omega_{3}\wedge\omega_{2},\quad d\omega_{2}=\omega_{1}\wedge\omega_{3},\quad d\omega_{3}=\omega_{1}\wedge\omega_{2},
\end{equation}
where $\omega_{i}=f_{i1}\,dx+f_{i2}\,dt$, $1\leq i\leq3$, are $1$-forms
with $f_{ij}$ smooth functions of $t$, $x$, $z$ and finite derivatives of $z$ with respect to $x$ and $t$, and such that $\omega_{1}\wedge\omega_{2}\neq0$
for generic solutions. In view of (\ref{struct}), for any generic
solution $z$ of $\mathcal{E}$ defined on an open subset
$U\subset\mathbb{R}^{2}$, the $1$-form $\omega_{3}=\omega_{12}$ defines the Levi-Civita connection of the corresponding pseudospherical
metric defined on $U$ by $I=(\omega_{1})^{2}+(\omega_{2})^{2}$.

The best known example of a \textbf{pss} equation is the sine-Gordon equation
\begin{equation}
z_{tx}=\sin z,\label{sG}
\end{equation}
with $1$-forms 
\begin{equation}
\omega_{1}=\frac{1}{\eta}\sin\left(z\right)\,dt,\quad\omega_{2}=\eta\,dx+\frac{1}{\eta}\cos\left(z\right)\,dt,\quad\omega_{3}=z_{x}\,dx,\label{sGcof2}
\end{equation}
where $\eta\in\mathbb{R}\setminus\{0\}$. In this case, one
has the metric given by $I=\omega_{1}^{2}+\omega_{2}^{2}=\frac{1}{\eta^{2}}dt^{2}+2\,\cos(z)\,dx\,dt+\eta^{2}dx^{2}$.
Other well-known differential equations that describe pseudospherical surfaces include for example, the KdV equation, the MKdV, the Burgers, the nonlinear Scr\"odinger, the short pulse equation, etc.

According to \cite{ChernTen}, \textbf{pss} equations have special properties that play an important role in the study of nonlinear evolution equations. In fact, \textbf{pss} equations  may have B\"acklund transformations, infinite hierarchies of conservation laws and they may also be solved by
applying the inverse scattering method (see for instance \cite{ChernTen},  \cite{BealsRabTen}
and references therein). Any {\bf pss} equation is the integrability condition of  a $2\times 2$ and  also  a $3\times 3$ linear problem determined by the 1-forms $\omega_i$, $i=1,2,3$ satisfying \eqref{struct}. In fact, the linear problem is given by 
\begin{equation}\label{linearproblem}
\frac{\partial V}{\partial x}=X\, V, \qquad \frac{\partial V}{\partial t}=T\, V,
\end{equation}
where
\[
X=\frac{1}{2}\left(\begin{array}{cc}
f_{21} & f_{11}-f_{31}\\
f_{11} +f_{31}& -f_{21}
\end{array}\right),\qquad
T=\frac{1}{2}\left(\begin{array}{cc}
f_{22} & f_{12}-f_{32}\\
f_{12} +f_{32}& -f_{22}
\end{array}\right)
\]
or 
\[
X=\left(\begin{array}{ccc}
0 & f_{11} & f_{21} \\
f_{11} & 0 & f_{31} \\
f_{21} & -f_{31} & 0
\end{array}\right),\qquad\qquad 
T=\left(\begin{array}{ccc}
0 & f_{12} & f_{22} \\
f_{12} & 0 & f_{32} \\
f_{22} & -f_{32} & 0
\end{array}\right).
\] 
 For recent developments on the study and
classification of \textbf{pss} equations, the reader is referred to \cite{Tarcisio-Niky,
Tarcisio-Keti, Catalano-Tenenblat, Catalano-Tarcisio-Keti,Catalano-Silva,
Diego-Luis,KT, 
NNK1,NNK2,NNK3,
RabTen2}.

Any pseudospherical surface described by a \textbf{pss} equation
$\mathcal{E}$ admits a local isometric immersion into $\mathbb{E}^{3}$
: this  is due to the classical result that a pseudospherical surface
always admits a local isometric immersion into $\mathbb{E}^{3}$ .
Hence, in view of Bonnet theorem, to any generic solution $z$
of $\mathcal{E}$ it is associated a pair $\left(I[z],II[z]\right)$
of first and second fundamental forms, which satisfy the Gauss--Codazzi
equations and describe a local isometric immersion into $\mathbb{E}^{3}$
of the associated pseudospherical surface. Recall that (see \cite{NNK2}
and also \cite{Diego-Luis,NNK1,NNK3,Tarcisio-Niky}) by introducing
the $1$-forms 
\begin{equation}\label{om1323}
\omega_{13}=a\omega_{1}+b\omega_{2},\qquad\omega_{23}=b\omega_{1}+c\omega_{2},
\end{equation}
the {\em Gauss--Codazzi equations}  read 
\[
ac-b^{2}=-1,\qquad d\omega_{13}=\omega_{12}\wedge\omega_{23},\qquad d\omega_{23}=\omega_{21}\wedge\omega_{13}
\]
and the second fundamental form writes as 
\[
II=\omega_{1}\cdot\omega_{13}+\omega_{2}\cdot\omega_{23}.
\]
In particular, $H=(a+c)/2$ gives the mean curvature of the isometric
immersion.

\selectlanguage{english}%
Thus for any generic solution $z$ of $\mathcal{E}$, defined on some
open subset $U\subset\mathbb{R}^{2}$, there always exists a
triple $\{a,b,c\}$ of differentiable functions that is locally defined
on $U$ and satisfies the Gauss--Codazzi equations for the pseudospherical
surface described by $z$. However, in general the dependence of $\{a,b,c\}$, and hence
of  $\left(I,II\right)$, on $z$ may be quite complicated, and it is
not guaranteed that it depends on a jet of finite order of $z$, i.e.,  on $x$, $t$ and a finite 
number of derivatives of $z$ with respect to $x$ and $t$. 

\selectlanguage{british}%

An example of a \textbf{pss} equation which admits local isometric immersions with $\{ a,b,c\}$ (and hence principal curvatures) depending on a jet of finite order of $z$, is provided by the sine-Gordon equation: indeed by \foreignlanguage{english}{taking }$\omega_{1}$,
$\omega_{2}$ and $\omega_{3}$ as in (\ref{sGcof2}), one has 
the first and second fundamental forms given by 
$$
I=\frac{1}{\eta^{2}}dt^{2}+2\,\cos(z)\,dx\,dt+\eta^{2}dx^{2} \qquad 
 II=\pm2\,\sin(z)\,dx\,dt, 
 $$
satisfying the  Gauss--Codazzi equations. The corresponding 1-forms $\omega_{13}$ and $\omega_{23}$ are given by \eqref{om1323}, where 
$$a=\pm2/\tan(z), \qquad 
b=\mp1, \qquad c=0. 
$$

In  a series of recent papers \cite{Tarcisio-Niky,Diego-Luis,NNK1,NNK2,NNK3},
the existence of such local isometric immersions has been
investigated for the families of \textbf{pss} equations previously studied in the papers
\cite{Tarcisio-Keti,Catalano-Silva,ChernTen,KT,RabTen2}. Surprisingly,
those papers showed that, except for the sine-Gordon equation,  all local isometric immersions admitted by
\textbf{pss} equations described in these papers have triples, $\{a,b,c\}$ which are  ``universal''
in the sense that they only depend on $x$ and $t$. 
Hence, contrary
to the case of the sine-Gordon equation, for any \textbf{pss} equation $\mathcal{E}$
of the type described in \cite{Tarcisio-Keti,Catalano-Silva,ChernTen,KT,RabTen2} the principal curvatures of the immersion are functions of $x$, $t$ and do not depend on the generic solutions. This fact shows that the 
sine-Gordon equation occupies a particularly special place amongst all \textbf{pss} equations,
and clearly motivates the search for other examples of \textbf{pss} equations
which admit local isometric immersions whose principal curvatures
 depend explicitly on their generic solutions. Families of equations of this type are provided in Section 2.

Motivated by our recent paper \cite{Catalano-Tarcisio-Keti}, we start by introducing a family of \textbf{pss} differential equations depending  on four arbitrary differentiable functions which satisfy generic conditions (see Proposition \ref{Prop1}). Then, by specializing these functions, we obtain a family of \textbf{pss} equations that admit local isometric immersions with  principal curvatures depending on the first order jet of solutions of the \textbf{pss} equations (see Corollary \ref{Cor1} and Proposition \ref{Prop2}). Next, we provide explicit examples, which include a generalization of the short pulse equation. 

We point out that, since the nineteenth century, the sine-Gordon equation has been classically related to surfaces of constant negative curvature; moreover, it also appeared in a number of physical applications, including the relativistic field theory and  the mechanical transmission lines.  
The short pulse equation, on the other hand, first appeared in the study of {\bf pss} equations \cite{BealsRabTen,Rab}, and only recently it has found physical applications in the description of ultra-short light pulses in silica optical fibers \cite{SchaferWayne}. It would be interesting to find out if the families of equations given in Section \ref{sec2}, contain other  equations, besides the short pulse equation, also describe physical phenomena. 

\section{Equations}

\label{sec2}

Recently, in \cite{Catalano-Tarcisio-Keti} we obtained a classification
for the class of quasilinear second order partial differential equations
\begin{equation*}
z_{tt}=A(z,z_{x},z_{t})\,z_{xx}+B(z,z_{x},z_{t})\,z_{xt}+C(z,z_{x},z_{t}), 
\end{equation*}
which describes either pseudospherical surfaces or spherical surfaces,
under suitable conditions for the functions $f_{ij}$.

The following result is a consequence of Theorem 3.6 from \cite{Catalano-Tarcisio-Keti}.

\begin{prop}\label{Prop1} A second order partial differential equation
\begin{equation}\label{eq:basic equation}
z_{xt}= \frac{1}{f_{12,z_t}-f_{11,z_x}}(f_{11,z_t}\,z_{tt}-f_{12,z_x}\,z_{xx}+  f_{11,z}\,z_t - f_{12,z}\,z_x + f_{31}f_{22}-f_{21}f_{32}), 
\end{equation}
where $f_{ij}=f_{ij}(z,z_{t},z_{x})$, $1\leq i \leq 3$ and $1\leq j \leq 2$, are differentiable functions such that $f_{12,z_t}-f_{11,z_x}\neq 0$ and $f_{21,z_{t}}=f_{31,z_{t}}=0$,  describes \textbf{pss}  if, and only if, 
\begin{equation}\label{f11f12}
\left(\begin{array}{ccc}
f_{11}\\
\\
f_{12}
\end{array}\right)=\dfrac{1}{\Delta}\left(\begin{array}{ccc}
-\phi z_{x}-\psi_{21} &  & \varphi z_{x}+\psi_{31}\\
\\
-\phi z_{t}-\psi_{22} &  & \varphi z_{t}+\psi_{32}
\end{array}\right)\left(\begin{array}{ccc}
\psi_{21,z}\qquad-\psi_{22,z}\\
\\
\psi_{31,z}\qquad-\psi_{32,z}
\end{array}\right)\left(\begin{array}{ccc}
z_{t}\\
\\
z_{x}
\end{array}\right)
\end{equation}
and
\begin{equation}\label{f21etc}
\begin{array}{l}
f_{21}=\phi z_{x}+\psi_{21},\vspace{5pt}\\
f_{31}=\varphi z_{x}+\psi_{31},
\end{array}\qquad\begin{array}{l}
f_{22}=\phi z_{t}+\psi_{22},\vspace{5pt}\\
f_{32}=\varphi z_{t}+\psi_{32},
\end{array}
\end{equation}
with $f_{11,z_{t}}\neq 0$, $\Delta=f_{32}f_{21}-f_{31}f_{22}\neq 0$ and $f_{21}f_{12}-f_{11}f_{22}\neq 0$, where $\phi$, $\varphi$ and $\psi_{rs}$, $r=2,3$, $s=1,2$ are differentiable functions of $z$.
\end{prop}

\begin{rem}
Since $\Delta = f_{32}f_{21} - f_{31}f_{22}$, it is a straightforward computation to check that \eqref{eq:basic equation} can be written in the more compact form
\begin{eqnarray*}\label{compact_eq}
D_x(f_{12}) - D_t(f_{11}) + \Delta = 0,
\end{eqnarray*}
where $D_x$ and $D_t$ are total derivatives.
\end{rem}

The proposition above provides a huge family of \textbf{pss} differential equations. In fact, they are obtained by choosing six arbitrary differentiable functions of $z$, satisfying generic conditions, namely $\phi$, $\varphi$ and $\psi_{rs}$, $r=2,3$, $s=1,2$, which determine  $f_{ij}$ by \eqref{f11f12} and \eqref{f21etc}. \\

In what follows, we will specialize these functions to
{\small $$
\phi(z) = \dfrac{2(F_1(z)+\lambda^2\,F_2(z))}{\lambda\,z^2},\quad \varphi(z) = -2\epsilon\,\dfrac{F_1(z)}{z},\quad \psi_{21}(z) = \lambda, \quad \psi_{22}(z) = \lambda\,\dfrac{z^2}{2}+\dfrac{1}{\lambda}, \quad \psi_{31}(z) = 0, \quad \psi_{32}(z) = -\epsilon\,z,
$$}\noindent
where $\epsilon=\pm 1$; $F_1(z)$, $F_2(z)$ are differentiable functions and $\lambda$ $\in$ $\mathbb{R}\setminus \{0\}$ is a free parameter. Then, as a consequence of Proposition \ref{Prop1}, we get the following

\begin{cor}\label{Cor1}{\em 
The differential equation }
\begin{equation}\label{4a}
D_x\left[\dfrac{(z^4+4F_2\,z_t)\,z_x}{2 G}\right] - D_t\left[\dfrac{(z^2+2F_2\,z_x)\,z_x}
{G}\right]+\dfrac{1}{z}G=0,
\end{equation}
where 
\[
G(z,z_x,z_t)= (F_1\,z^2-2F_2)\,z_x-2F_1\,z_t-z^2, 
\]
{\em  $F_1(z)$ and $F_2(z)$ are differentiable functions and $G\neq 0$, describes \textbf{pss} equations with
\begin{equation}\label{fijcor1}
\begin{array}{lll}
&&\omega_1 = \dfrac{\epsilon\lambda\,z_x}{2G}\left[2(z^2 + 2F_2\,z_x)\, dx + (z^4 + 4F_2\,z_t)\,dt\right], \vspace*{0.4cm}\\
&&\omega_2 = \left[\dfrac{2(\lambda^2 F_2 + F1)}{\lambda z^2}\,z_x+\lambda \right]\,dx + \left[\dfrac{2(\lambda^2 F_2 + F1)}{\lambda z^2}\,z_t + \dfrac{\lambda}{2}\,z^2 + \dfrac{1}{\lambda}\right]\,dt,\vspace*{0.4cm}\\ 
&&\omega_3 = -2\epsilon F_1\,\dfrac{z_x}{z}\,dx + \left(-2\epsilon F_1\,\dfrac{z_t}{z}-\epsilon z\right)\,dt,
\end{array}
\end{equation}
and  $\lambda\in \mathbb{R}\setminus \left\lbrace 0 \right\rbrace$.} 
\end{cor}


The next result shows that \eqref{4a} is a \textbf{pss} equation that admits local isometric immersions, whose principal curvatures depend on a finite order jet of $z$. 

\begin{prop}\label{Prop2}
For any solution $z$ of \eqref{4a} which  satisfies $z_x\neq 0$, the pseudospherical metric  $I=\omega_{1}^{2}+\omega_{2}^{2}$ has an isometric immersion in $\mathbb{E}^3$,  with $a$, $b$ and $c$ depending on the first order jet of $z$,  which is  given by
\[
\omega_{13} =\frac{ -2\epsilon}{z²z_x}\,G\,\omega_1+\omega_2, \qquad \omega_{23} = \omega_1.
\]
\end{prop}

\noindent
\textbf{Proof}. Recalling that $\omega_3=\omega_{12}$ is the connection form of $I=\omega_{1}^{2}+\omega_{2}^{2}$, a straightforward computation with $\omega_1$, $\omega_2$ and $\omega_3$ as in Corollary \ref{Cor1} entails that
$$
d\omega_{13} = \omega_{12}\wedge\omega_{23}, \quad d\omega_{23} = \omega_{21}\wedge\omega_{13}\quad \textnormal{and} \quad ac-b^2 = -1.
$$

\hfill
$\square$

Before providing  particular examples, we observe that 
 for any choice of $F_1(z)$ and $F_2(z)$ as in Corollary \ref{Cor1}, the  corresponding differential  equation \eqref{4a} 
  is the compatibility condition of a linear problem (Lax pair or zero curvature condition) described  by \eqref{linearproblem},  where the  functions $f_{ij}$ are explicitly given by  \eqref{fijcor1}, in terms of $z,\, z_x,\, z_t$ and a non zero parameter $\lambda$, since  $\omega_i=f_{i1}\,dx+f_{i2}\,dt$.     
  
Moreover, one may obtain infinite conservation laws for \eqref{4a} as a consequence of  the following geometric properties. Assume that a differential equation $\cal{E}$ for $z(x,t)$ describes   pseudospherical surfaces, with associated 1-forms $\omega_1,\, \omega_2,\, \omega_3$, then $z$ is a solution of $\cal{E}$ if and only if $\omega_3-d\rho\,\omega_1+\sin\rho\,\omega_1+\cos\rho\, \omega_2=0$ 
is completely integrable for $\rho$.  For each solution $z$ of $\cal{E}$ and corresponding solution $\rho$, the 1-form $\cos\rho\, \omega_1 -\sin\rho\, \omega_2$ is closed (see the  proof of Proposition 4.2 in \cite{ChernTen}). Whenever the 1-forms  $\omega_i$, $i=1,2,3$ are analytic on a parameter $\lambda$, then $\rho$  is also analytic on the parameter and hence the conservation laws may be obtained by the closed form written as a series in $\lambda$.    

\vspace{.1in}

We now exhibit   a class of examples of \eqref{4a} which includes the short pulse equation.  

\begin{example}\label{exe2.1a} 
In view of \eqref{4a}, by considering $F_1(z) = \ell(z)\,z^2$ and $F_2(z) = 0$, where $\ell(z)$ is a differentiable function, one has  
\begin{eqnarray*}
z_{xt} = 
\frac{1}{2[\ell\,(z^2z_x +2z_t)+1 ]}\left\{4\ell z_x z_{tt}+z^2\left(1+2\ell z_t\right)z_{xx}- \ell' (4z^2z_xz_t-z^4z_x^2-4z_t^2) z_x + \right. \\
\left. + 2z z_x^2  
-2z \left[\ell (z^2z_x-2z_t)-1\right]^3
 \right\}. 
\end{eqnarray*}

\noindent
This is a {\bf pss} differential equation, whose generic solutions define pseudospherical metrics admitting local isometric immersions in $\mathbb{R}^3$, with first and second fundamental forms provided by
\begin{equation} \label{omexe}
\begin{array}{lll}
&&\omega_1 = \dfrac{\epsilon\lambda\,z_x}{\ell\,(z^2z_x-2\,z_t)-1}
\left( dx+ \dfrac{z^2}{2}\, dt \right),\\
,\vspace*{0.4cm}\\
&&\omega_2 = \left(\dfrac{2\ell}{\lambda}\,z_x+\lambda \right)\,dx + \left(\dfrac{2\ell}{\lambda}\,z_t + \dfrac{\lambda}{2}\,z^2 + \dfrac{1}{\lambda}\right)\,dt,\vspace*{0.4cm}\\ 
&&\omega_3 = -2\epsilon\, \ell\,z \,z_x\,dx -\epsilon\,z 
\left(2 \ell\,z_t+1 \right)\,dt,
\end{array}
\end{equation}
and
\begin{eqnarray*}
\omega_{13} = -\dfrac{2\epsilon}{z_x}\left[ \ell(z^2z_x-2z_t)-1 \right]\,\omega_1 + \omega_2, \quad \omega_{23} = \omega_1.
\end{eqnarray*}
\end{example}
\begin{example}\label{exe2.1} 
By choosing $\ell(z) = 0$ in the previous example (i.e. $F_1=F_2=0$ in \eqref{4a}), one obtains the \textit{short pulse equation}
$$
z_{xt} = z + \frac{1}{6}\,(z^3)_{,xx}, 
$$
also known as the Rabelo's cubic equation \cite{Rab}, with first and second fundamental
forms provided by
\[
\omega_1 = -\epsilon\lambda\,z_x\, dx - \dfrac{\epsilon\lambda}{2}\, z^2\,z_x\,dt,\qquad
\omega_2 = \lambda\,dx + \left(\dfrac{\lambda}{2}\,z^2 + \dfrac{1}{\lambda}\right)\,dt,\qquad 
\omega_3 = -\epsilon z\,dt,
\]
and
\begin{eqnarray*}
\omega_{13} = \dfrac{2\epsilon}{z_x}\,\omega_1 + \omega_2, \quad \omega_{23} = \omega_1.
\end{eqnarray*}
\end{example}
\begin{example}\label{exe2.2} 
By choosing $\ell(z)=m\in\mathbb{R}$ in Example \ref{exe2.1a}, one obtains the following equation
 \begin{eqnarray*}
z_{xt} = 
\frac{1}{2[m\,(z^2z_x +2z_t)+1 ]}\left\{4m z_x z_{tt}+z^2\left(1+2mz_t\right)z_{xx}
+ 2z z_x^2  -2z \left[m (z^2z_x-2z_t)-1\right]^3
 \right\},  
\end{eqnarray*}
 which is a {\bf pss} equation with associated forms 
 as in  \eqref{omexe}, where $\ell=m$. 
\end{example}

\noindent
Thus, the family of equations given by Example \ref{exe2.1a} (in particular by Example \ref{exe2.2}) provides {\bf pss} equations that  generalize the short pulse equation and admit local isometric immersions in $\mathbb{R}^3$ whose principal curvatures depend on the first order jet of its solutions. We note here that other generalizations of the short pulse equation have been discussed and investigated in \cite{Catalano-Tarcisio-Keti, HNW,Sak1}.

 
 \section*{Acknowledgements}
We would like to thank the referees for their valuable comments that contributed to the improvement of our paper. 

\vspace{.1in}  
 
DCF was partially supported by CNPq, grant 310577/2015-2 and grant 422906/2016-6. TCS was partially supported by FAPDF/Brazil, grant 0193.001346/2016, CAPES/Brazil-Finance Code 001 and by CNPq, grant 422906/2016-6. KT was partially supported by CNPq, grant 312462/2014-0, CAPES/Brazil-Finance Code 001 and FAPDF /Brazil, grant 0193.001346/2016.


\begin{thebibliography}{10}
\bibitem{BealsRabTen} R. Beals, M.Rabelo, K.Tenenblat, Backl\"und transformations and inverse 
scattering solutions for some pseudospherical surface equations, Stud. Appl. Math. \textbf{81}, No. 2
(1989)  125--151. 

\bibitem{Tarcisio-Niky} T. Castro Silva and N. Kamran, Third order
differential equations and local isometric immersions of pseudospherical
surfaces, Communications in Contemporary Math. \textbf{18}, No. 6
(2016) 1650021 (41 pages).

\bibitem{Tarcisio-Keti} T. Castro Silva and K. Tenenblat, Third order
differential equations describing pseudospherical surfaces, J. Differential
Equations \textbf{259} (2015) 4897--4923.

\bibitem{Catalano-Tenenblat}D. Catalano Ferraioli and K. Tenenblat,
Fourth order evolution equations which describe pseudospherical surfaces,
J. Differential Equations \textbf{257} (2014) 3165--3199.

\bibitem{Catalano-Tarcisio-Keti} D. Catalano Ferraioli, T. Castro
Silva and K. Tenenblat, A class of quasilinear second order partial
differential equations which describe spherical or pseudospherical
surfaces, J. Differential Equations \textbf{268} (2020) 7164--7182.

\bibitem{Catalano-Silva}D. Catalano Ferraioli, L. A. de Oliveira
Silva, Second order evolution equations which describe pseudospherical
surfaces, J. Differential Equations {\bf 260} (2016) 8072--8108.

\bibitem{Diego-Luis} D. Catalano Ferraioli and L. A. de Oliveira,
Local isometric immersions of pseudospherical surfaces described by
evolution equations in conservation law form, J.Math. Anal. Appl. \textbf{446} (2017)
1606--1631.

\bibitem{ChernTen}S. S. Chern, K. Tenenblat, Pseudo-spherical surfaces
and evolution equations, \textit{Stud. Appl. Math.} \textbf{74}, (1986)
55--83.

\bibitem{HNW} A. N. W. Hone, V. Novikov and J. P. Wang, Generalizations of the short pulse equation, \textit{Lett. Math. Phys.} \textbf{108}, (2018)
927--947.

\bibitem{KT}N. Kamran and K. Tenenblat, On differential equations
describing pseudospherical surfaces, \textit{J. Differential Equations}
\textbf{115} (1995) 75--98.

\bibitem{NNK1} N. Kahouadji, N. Kamran and K. Tenenblat, Local isometric
immersions of pseudo-spherical surfaces and evolution equations, Fields
Inst. Commun. \textbf{75} (2015) 369--381.

\bibitem{NNK2} N. Kahouadji, N. Kamran and K. Tenenblat, Second-order
equations and local isometric immersions of pseudo-spherical surfaces.
Comm. Anal. Geom. \textbf{24} (3) (2016) 605--643.

\bibitem{NNK3} N. Kahouadji, N. Kamran and K. Tenenblat, Local isometric
immersions of pseudo-spherical surfaces and kth-order evolution equations,
Communications in Contemporary Math. \textbf{21}, No. 4 (2019) 1850025
(21 pages).

\bibitem{Rab} M. Rabelo, On equations which describe pseudospherical surfaces ,
Stud. Appl. Math. \textbf{81} (1989) 221--248.


\bibitem{RabTen2} M. Rabelo, K. Tenenblat, A classification of pseudospherical
surface equations of type $u_{t}=u_{xxx}+G(u,u_{x},u_{xx})$,
J. Math. Phys. \textbf{33} (1992) 537--549.

\bibitem{Sak1} A. Sakovich, S. Sakovich, The Short Pulse Equation
is Integrable, J. Phys. Soc. Jpn. \textbf{74},(2005) 239--241.

\bibitem{SchaferWayne} T. Sch\"afer, C.E. Wayne, Propagation of ultra short optical pulses in cubic nonlinear media,  Physica D \textbf{196} (2004) 90--105. 

\end{thebibliography}
\end{document}